\documentclass{article}
\usepackage{authblk}
\usepackage{biblatex}

\usepackage{booktabs}

\usepackage[]{blindtext}
\usepackage{amsmath}
\usepackage{amssymb}

\title{A Random-Player Game and Derangement Numbers}
\author{Yehonatan Fridman}
\affil{Ben-Gurion University, NRCN, Israel \thanks{fridyeh@post.bgu.ac.il}}
\date{} 

\begin{document}

\maketitle

\begin{abstract}
 Consider the following game between a random player $\mathfrak{R}$ and a deterministic player $\mathfrak{D}$. There is a pile of $n$ elements at the beginning. The rules for playing are as follows: In each turn of $\mathfrak{R}$, if the pile contains exactly $m$ elements, $\mathfrak{R}$ removes $k$ elements from the pile, where $k$ is independently identically distributed from \{$1, \ldots, m\}$. In each turn of $\mathfrak{D}$, $\mathfrak{D}$ removes only one element. The winner is the player that, at the end of its round, has no elements remaining. $\mathfrak{R}$ starts first to play. This short paper shows that $\lim_{n\to\infty} [D_n:=p(\mathfrak{D} \text{ wins when game initialized with } n \text{ elements})] = e^{-1}$; and more specifically, $D_n = \frac{d_n}{n!}$, where $d_n$ is the $n$-th derangement number.
\end{abstract}

\subsection*{keywords}
Combinatorics, Random Games, Derangement Numbers

\section{Game Definition}
This paper introduces a game between a random player $\mathfrak{R}$ and a deterministic player $\mathfrak{D}$, where a pile of $n$ elements initiates the game. The players take turns alternately, with player $\mathfrak{R}$ starting the game. During $\mathfrak{R}$'s turn, it removes $k$ elements from the pile, where $k$ is independently identically distributed (i.i.d.) from $\{1, \ldots, m\}$, given that the pile contains exactly $m$ elements at that moment. Conversely, player $\mathfrak{D}$ removes only one element during its turn (and hence is deterministic). The game concludes when one player clears the pile completely, thus winning. Let denote the probability of $\mathfrak{D}$ winning, when the game is initialized with $n$ elements, as $D_n$, and the probability of $\mathfrak{R}$ winning as $R_n$. It is clear that $D_n + R_n = 1$. The subsequent section introduces a recursive formula to calculate $R_n$, followed by its solution, thereby providing an explicit method to determine $R_n$ (and consequently $D_n$).

\section{Calculating the Probability of $\mathfrak{R}$ Winning}
It is clear to see that $R_1=1$ and $R_2=0.5$. A recursive relation to describe $R_n$ (for $n\geq 3$) can be achieved by conditioning on all the options of the first turn of $\mathfrak{R}$. Each event of removing $k \in \{1, \ldots, n\}$ elements in the first turn of $\mathfrak{R}$ occurs with probability $1/n$, and winning depends recursively on the game's progress when $\mathfrak{R}$ is given to play with $n-k-1$ elements remaining. Special cases are for $k=n$ and $k=n-1$, which determine the game's result after the first turn of $\mathfrak{R}$: if $k=n$, then $\mathfrak{R}$ wins, and if $k=n-1$, then $\mathfrak{R}$ loses. Therefore, 
\begin{equation}
\label{eq1}
R_n = \frac{1}{n} \sum_{k=1}^{n-2} R_{n-k-1} + \frac{1}{n}\cdot0 + \frac{1}{n}\cdot1 = \frac{1}{n} \sum_{k=1}^{n-2} R_k + \frac{1}{n} = \frac{1}{n}\cdot\Biggl(\sum_{k=0}^{n-2} R_k + 1\Biggl)\tag{1}
\end{equation}

when the last equation holds when defining $R_0=0$ (thus the random player is said to be losing if the pile is initialized with no elements). Further development of Equation (\ref{eq1}) yields:

\begin{equation}
\label{eq2}
R_n = \frac{1}{n}\cdot\Biggl(\sum_{k=0}^{n-2} R_k + 1\Biggl) = \frac{1}{n}\Biggl(1+ R_{n-2}+\sum_{k=0}^{n-3} R_k\Biggl)=\frac{1}{n}\Biggl( R_{n-2} + (n-1)R_{n-1}\Biggl)\tag{2}
\end{equation}

The solution of this recursive formula can be investigated using generating functions. Let define $\mathcal{R}(x)$ be the generating function of $R_n$. By definition:

\begin{equation}
\label{eq3}
    \mathcal{R}(x) = \sum_{n=0}^{\infty}R_nx^n \Longrightarrow{} \mathcal{R}'(x) = \sum_{n=1}^{\infty}n R_nx^{n-1} \Longrightarrow{} x\mathcal{R}'(x) = \sum_{n=0}^{\infty}n R_nx^n \tag{3}
\end{equation}

Hence, $x\mathcal{R}'(x)$ is the generating function of $nR_n$. From Equation (\ref{eq2}) we get:

\begin{equation}
\label{eq4}
\begin{split}
    nR_n &= R_{n-2} + (n-1)R_{n-1} \Longrightarrow \\
     & \sum_{n=2}^{\infty}nR_nx^n = \sum_{n=2}^{\infty}R_{n-2}x^n + \sum_{n=2}^{\infty}(n-1)R_{n-1}x^n \Longrightarrow \\ & \sum_{n=0}^{\infty}nR_nx^n - R_1x = x^2\sum_{n=2}^{\infty}R_{n-2}x^{n-2} + x\sum_{n=2}^{\infty}(n-1)R_{n-1}x^{n-1} = \\
     & x^2\sum_{n=0}^{\infty}R_nx^n + x\sum_{n=1}^{\infty}nR_nx^n = x^2\sum_{n=0}^{\infty}R_nx^n + x\sum_{n=0}^{\infty}nR_nx^n
     \end{split} \tag{4}
\end{equation}

Now, using the colerations in Equation (\ref{eq3}) for $\mathcal{R}(x)$ and $x\mathcal{R}'(x)$, the following holds directly from Equation (\ref{eq4}):

\begin{equation}
\label{eq5}
\begin{split}
    x\mathcal{R}'(x) &-x = x^2\mathcal{R}(x)+x^2\mathcal{R}'(x) \Longrightarrow \mathcal{R}'(x) -1 = x\mathcal{R}(x)+x\mathcal{R}'(x) \Longrightarrow \\
    & (x-1)\mathcal{R}'(x) = -1 -x\mathcal{R}(x) \Longrightarrow \mathcal{R}'(x) + \frac{x\mathcal{R}(x)}{x-1}= \frac{-1}{x-1}
    \end{split} \tag{5}
\end{equation}

To solve this differential equation, one can choose the following integrating factor: 
\begin{equation*}
exp \Biggl(\int \frac{x}{x-1} \Biggl) = e^x(x-1)
\end{equation*}

Which leads to the equivalent equation:

\begin{equation}
\begin{split}
e^x(x-1)&\mathcal{R}'(x) + xe^x\mathcal{R}(x) = -e^x \Longrightarrow \frac{d}{dx}\Bigl( e^x(x-1)\mathcal{R}(x)\Bigr)=-e^x \Longrightarrow \notag \\
& e^x(x-1)\mathcal{R}(x) = -e^x+C \Longrightarrow \mathcal{R}(x)= \frac{Ce^{-x}}{x-1} -\frac{1}{x-1} 
\end{split}
\label{eq6} \tag{6}
\end{equation}

and because $R_0=0$ (as defined before), we get $C=1$. Therefore:

\begin{equation}
\mathcal{R}(x) = \frac{e^{-x}-1}{x-1} \Longrightarrow \mathcal{R}(x)=\mathcal{F}(x)\mathcal{G}(x) \tag{7}
\label{eq7}
\end{equation}

where $\mathcal{F}(x):=\frac{1}{x-1}$ and $\mathcal{G}(x):=e^{-x}-1$.

\vspace{5mm}

It is known that:
\begin{equation}
\mathcal{F}(x) = -\sum_{n=0}^{\infty}x^n \quad \text{and} \quad \mathcal{G}(x) = \sum_{n=0}^{\infty}\frac{(-x)^n}{n!}-1 \tag{8}
\label{eq8}
\end{equation}

From Equations (\ref{eq7}) and (\ref{eq8}) we get:

\begin{equation}
\begin{split}
\mathcal{R}(x) &= -\Biggl(\sum_{n=0}^{\infty}x^n \Biggl)\Biggl(\sum_{n=0}^{\infty}\frac{(-x)^n}{n!}-1\Bigg)
= \Biggl(\sum_{n=0}^{\infty}x^n \Biggl)\Biggl(1-\sum_{n=0}^{\infty}\frac{(-1)^nx^n}{n!}\Bigg) \\
&= \sum_{k=0}^{\infty}\Biggl(1-\sum^{k}_{l=0}\frac{(-1)^l}{l!}\Biggl)x^k
\end{split}
\label{eq9} \tag{9} 
\end{equation}

And because $\mathcal{R}(x)$ is the generative function of $R_n$, we finally get: 
\begin{equation}
\begin{split}
R_n = 1-\sum^{n}_{k=0}\frac{(-1)^k}{k!} \quad  \text{ and } \quad D_n = 1-R_n = \sum^{n}_{k=0}\frac{(-1)^k}{k!}
\end{split}
\label{eq10} \tag{10} 
\end{equation}
It is a well-known result that $\frac{d_n}{n!} = \sum_{k=0}^{n}\frac{(-1)^k}{k!}$, where $d_n$ represents the n-th derangement number~\cite{carlitz1978number}. Therefore, we conclude that $D_n = \frac{d_n}{n!}$ for every $n$ (one can verify that this also holds for the base cases $n=1,2$). In other words, the probability of $\mathfrak{D}$ winning a game with $n$ elements is equivalent to the probability of obtaining a fixed point-free permutation when randomly selecting a permutation of order $n$.
This result further explains the behavior of the described game as the pile size approaches infinity: $\lim_{n\to\infty} [D_n := p(\mathfrak{D} \text{ wins when the game is initialized with } n \text{ elements})] = e^{-1}$.

\subsection{Alternative Solution} \label{alter}
Another solution for the inductive formula may be considered simpler and more elegant, and it is provided by the following derivation. We start again with:

\begin{equation}
    R_n=\frac{1}{n}\Bigl(  R_{n-2}+(n-1)R_{n-1}\Bigl)
    \label{eq11} \tag{11}
\end{equation}

Rearranging gives:

\begin{equation}
    R_n-R_{n-1}=\frac{1}{n}\Bigl( R_{n-2}-R_{n-1}\Bigl)
 \Longrightarrow a_n=\frac{-1}{n}a_{n-1}
 \label{eq12} \tag{12}
\end{equation}

 where $a_n:=R_n-R_{n-1}$ and $a_1:=1$ (remember that $R_0=0$ and $R_1=1$).

\vspace{5mm}
 
 The solution is clearly 

 \begin{equation}
    a_n=\frac{(-1)^{n+1}}{n!}
    \label{eq13} \tag{13}
\end{equation}

Therefore,

 \begin{equation}
   \begin{split}
    R_n &= R_0+a_1+a_2+ \dots + a_n \\
    &= 0 + 1 -\frac{1}{2!} + \frac{1}{3!} + \dots + \frac{(-1)^{n+1}}{n!} = 1- \sum^{n}_{k=0}\frac{(-1)^k}{k!}
\end{split}
    \label{eq14} \tag{14}
\end{equation}

\section{Conclusion}
This paper presents a straightforward random-player game, where the probability of the random player losing (or the deterministic player winning) is equivalent to $\frac{d_n}{n!}$, when ${d_n}$ is the n-th derangement number~\cite{carlitz1978number}. As mentioned in the previous section, this is equivalent to stating that the probability of $\mathfrak{D}$ winning a game with $n$ elements corresponds to the likelihood of selecting a permutation without fixed points when randomly choosing a permutation of order $n$. The proof offered in this paper is analytical rather than combinatorial. It would be intriguing to explore a combinatorial explanation for the close association between this game and the well-investigated derangement combinatorial problem.

\subsection*{Acknowledgements}
The author would like to thank Kent E. Morrison\footnote{Emeritus Professor of Mathematics
Cal Poly, San Luis Obispo.} for suggesting the alternative solution listed in Section \ref{alter}.

\printbibliography

\end{document}